\documentclass{amsart}

\usepackage{amsthm,amssymb,latexsym,amsmath}
\usepackage{mathrsfs}
\usepackage{graphicx}
\usepackage{color}
\usepackage[latin1]{inputenc}

\def\C{\mathbb C}

\renewcommand{\dim}{\mathrm{dim}}

\newcommand{\OO}{\mathcal O}
\newcommand{\F}{\mathscr F}
\newcommand{\sH}{\mathscr H}
\newcommand{\fol}{\mathscr F}

\renewcommand{\O}{{\mathcal O}}

\newcommand{\PP}{\mathbb{P}}

\newcommand{\singf}{{\rm Sing}({\mathscr F})}

\newcommand{\x}{\times}

\newcommand{\Sing}{{\rm Sing}}
\newcommand{\cod}{{\rm cod}}

\newtheorem{lema}{Lemma}[section]
\newtheorem{cor}[lema]{Corollary}

\newtheorem*{teo1}{Theorem 1}

\newtheorem{prop}[lema]{Proposition}

\theoremstyle{definition}
\newtheorem{remark}[lema]{Remark}

\newenvironment{demostracion}
{\noindent {{\it  Proof.}}}

\begin{document}

\title{
Global residue formula for logarithmic indices of 
one-dimensional foliations
}


\begin{abstract}
We prove     a  global residual formula in terms of logarithmic indices  for  one-dimensional holomorphic foliations, with isolated singularities, and  logarithmic along   normal  crossings   divisors. We also give a formula  for the total sum of the logarithmic  indices   if    the singular set of the foliation is contained in  the invariant divisor.    As  an application,  we  provide a formula for the number of singularities in the complement of the invariant  divisor  on complex  projective spaces.  Finally,   we obtain  a Poincar\'e-Hopf type formula for singular normal projective  varieties.  
\end{abstract}

\author{Maur\'icio Corr\^ea}
\address{\noindent  Maur\'icio Corr\^ea,
\newline
\indent 
Universit\`a degli Studi di Bari, 
Dipartimento di Matematica,
Via E. Orabona 4, I-70125 Bari, Italy.
\newline
\indent 
Departamento de Matem\'atica \\
Universidade Federal de Minas Gerais\\
Av. Ant\^onio Carlos 6627 \\
30123-970, Belo Horizonte-MG, Brazil.
} \email{mauriciojr@ufmg.br, mauricio.barros@uniba.it}

\author{Diogo da Silva Machado}
\address{\noindent  Diogo da Silva Machado, 
\newline
\indent
Departamento de Matem\'atica \\
Universidade Federal de Vi\c cosa\\
Avenida Peter Henry Rolfs, s/n - Campus Universit\'ario \\
36570-900 Vi\c cosa- MG, Brazil.} \email{diogo.machado@ufv.br}

\subjclass{Primary 32S65, 37F75; secondary 14F05}

\keywords{Logarithmic foliations, Poincar\'e--Hopf type theorem, residues}


\maketitle

\hyphenation{des-cri-bed}



\section{Introduction}

A complex  $\partial$-manifold \cite{YN} is a complex manifold of the form 
$\tilde{X} = X -  D$, where $X$ is an $n$-dimensional complex compact manifold and  $ D\subset X$ is a  divisor which is called the  boundary divisor.
S. Iitaka in \cite{Lita}  proposed a version of Gauss-Bonnet theorem for  $\partial$-manifold \cite{YN}. Y. Norimatsu \cite{YN}, R. Silvotti \cite{RS} and P. Aluffi \cite{Aluffi} have proved independently    such Gauss-Bonnet type  theorem for the case  when the boundary divisor $D$ has normal crossings singularities. More precisely, they proved the following formula  
$$
\int_{X}c_n(T_X(-\log\,  D)) \,\,\, =\,\,\, \chi (\tilde{X}),
$$
 where $T_X(-\log\,  D)$ denotes the 
vector bundle  of logarithmic  vector fields along $ D$ and 
$
\chi (\tilde{X})
$
is the Euler characteristic of $\tilde{X}$.  
 In \cite{MD}   the authors have proved  the following  Poincar\'e-Hopf  type  index theorem 
$$
\int_{X}c_n(T_X(-\log\,  D)) \,\,\, =\,\,\, \chi (\tilde{X})= \sum_{x\in  \Sing(v) \cap  \tilde{X} } \mu_x(v),
$$
where   $v $ is a global holomorphic  vector field on $X$     tangent  to the  divisor $D$, 
$\mu_x(v)$ is  the Milnor number of $v$ at $x$ and  under the hypothesis that Aleksandrov's logarithmic indices of  $v$ along $D$  vanish. Let us recall  Aleksandrov's definition of logarithmic index   \cite{Ale_p} in the context of one-dimensional holomorphic foliations:

  Let $\fol$ be a one-dimensional holomorphic foliation, with isolated singularities,  on a    complex manifold    $X$  and logarithmic along a divisor  $D$, we refer to Section 2 for this notion.
  Let $\Sing(\fol)$ be the singular set of $\fol$.  Let  $v\in T_{X}(-\log \,D)|_{U}$ be  a germ of vector field on $(U,x)$ tangent to $\fol$. The interior multiplication $i_{v}$ induces the complex of logarithmic differential forms
$$
0 \longrightarrow \Omega^n_{X,x}(\log\, D) \stackrel{i_v}{\longrightarrow} \Omega^{n-1}_{X,x}(\log\, D)\stackrel{i_v}{\longrightarrow} \cdots\stackrel{i_v}{\longrightarrow} \Omega^{1}_{X,x}(\log\, D)\stackrel{i_v}{\longrightarrow} \OO_{n,x} \longrightarrow 0.
$$
Since all singularities of $v$ are isolated, the $i_v$-homology groups of the complex $\Omega^{\bullet}_{X,x}(\log\, D)$ are finite-dimensional vector spaces. Thus, the Euler characteristic
$$
\chi (\Omega^{\bullet}_{X,x}(\log\, D), i_v) = \displaystyle\sum_{i=0}^n(-1)^i \dim  H_i(\Omega^\bullet_{X,x} (\log\, D),i_v)
$$
\noindent of the complex of logarithmic differential forms is well defined.  Since this number does not depend on local representative $v$ of the foliation $\fol$ at $x$, we define the {\it logarithmic index} of $\fol$  at the point $x$ by
$$
Log(\fol,D,x):= \chi (\Omega^{\bullet}_{X,x}(\log\, D), i_v).
$$
This  index is  inspired by G\'omez-Mont's homological index \cite{a08}.

We prove the following Baum-Bott type residual formula \cite{BB1} in terms of the  logarithmic indices for one-dimensional holomorphic foliations, with isolated singularities, and logarithmic along normal  crossings divisors:

\begin{teo1} \label{teo1_1} 
Let $\fol$ be a one-dimensional foliation with isolated singularities and logarithmic along a normal crossings divisor $D$ on a complex compact manifold $X$.  Then
$$
\displaystyle\int_{X}c_{n}(T_{X}(-\log\, D)- T_{\fol}) =    \sum_{x\in  \Sing(\fol)\cap (X\setminus  D)} \mu_x(\fol) + \sum_{x\in \Sing(\fol) \cap  D }    Log(\fol,D,x),
$$
where $T_{\fol}$ denotes  the tangent bundle of $\fol$ and  $ \mu_x(\fol)$ is the Milnor number of  $\fol$ at $x$. 
\end{teo1}


We observe that Theorem \ref{teo1_1}  generalizes our previous result in \cite{MD}, where we have assumed that the  singularities of $\fol$ are non-degenerate. In \cite[Theorem 1, (ii)]{CFLM}
we have proved an analogous result  for global   holomorphic vector fields  on $X$ such that $D$  has isolated singularities.

We obtain the following consequence.

\begin{cor}  Let $\fol$ be a  one-dimensional    foliation on $X$, with isolated singularities and logarithmic along  a normal crossings  divisor $D$. 
\begin{enumerate}
\item[(i)] If 
$Log(\fol,D,x)=0$,  for all $x \in \Sing(\fol)\cap D$,
then
$$
\int_{X}c_{n}(T_{X}(-\log\, D)- T_{\fol}) = \sum_{x\in \Sing(\fol)\cap (X\setminus D)}\mu_x(\fol).
$$
\medskip
\item[(ii)]  If 
$ \Sing(\fol) \subset D$,  then
$$
\int_{X}c_{n}(T_{X}(-\log\, D)- T_{\fol}) = \sum_{x\in \singf \cap  D }  Log(\fol,D,x).
$$
\end{enumerate}

\end{cor}

Now, let $\fol$ be a one-dimensional  foliation on the projective space  $\mathbb{P}^n$.  As a consequence of the above formulas we can provide a formula for the number of singularities of the foliation $\fol$ in the complement of the invariant  divisor $D$. Moreover,  under the assumption
$ \Sing(\fol) \subset D$ we also have a formula in terms of the total sum of the logarithmic index along $D$.

\begin{cor}\label{cor-Pn}   Let $\fol$ be a   one-dimensional    foliation on  $\mathbb{P}^n$, of degree $d>1$, with isolated singularities and logarithmic along  a normal crossings  divisor $D=\displaystyle \cup_{i=1}^k D_i$. Denote $d_i:=\deg(D_i)$.
\begin{enumerate}
\item[(i)] If 
$Log(\fol,D,x)=0$,  for all $x \in \Sing(\fol)\cap D$,
then
$$
\sum_{x\in \Sing(\fol)\cap (X\setminus D)}\mu_x(\fol) = \displaystyle\sum^{n}_{i=0} {n+1\choose i}  \sigma_{n-i}(d_1,\dots,d_k, d-1), 
$$
\end{enumerate}
 where $\sigma_{n-i} $ is the     complete symmetric function of degree $(n-i)$ in the
variables $(d_1,\dots,d_k,d-1)$. 

\noindent  If, in addition,  $d=2$, $n\geq k$ and $d_i=1$ for all $i$,  then   $  \Sing(\fol) \subset D$. Otherwise, we have that   $ \Sing(\fol)\cap (X\setminus D)  \neq \emptyset .$ 
\begin{enumerate}
\item[(ii)]  If 
$ \Sing(\fol) \subset D$,  then
$$
 \sum_{x\in \singf \cap  D }  Log(\fol,D,x)= \displaystyle\sum^{n}_{i=0} {n+1\choose i}  \sigma_{n-i}(d_1,\dots,d_k, d-1). 
 $$
\end{enumerate}
 \noindent   If, in addition,  $d\neq 2$  then  $\fol$ has at least one degenerate  singular point.

\end{cor}

Cukierman,   Soares,   and Vainsencher in  \cite{CSV} have provided a  formula prescribing the number of isolated singularities of a codimension one logarithmic foliation on $\mathbb{P}^n$. We notice that the Corollary \ref{cor-Pn}  can be seen as a version  for one-dimensional foliations of  this result which provides    the same number of singularities  for  $n=2$ and if   the logarithmic  indices  along the divisor  vanish. 
 Also,  the part (i) of  Corollary \ref{cor-Pn} generalizes and improves    \cite[Theorem 3]{MD}.
    
Another interesting application of Theorem \ref{teo1_1} is a Poincar\'e-Hopf type formula for singular varieties.  As usual  we denote by  $T_Y:=\sH om( \Omega_Y^1, \O_Y)$ the tangent sheaf of $Y$, where $\Omega_Y^1$ is the sheaf of  K\"ahler differentials of $Y$. We obtain the following result.
\begin{cor}\label{singlogres}  
Let $Y$ be  a normal variety and  $v\in H^0(Y, T_Y)$  a holomorphic vector field with isolated singularities such that $\Sing(v)\cap \Sing(Y)=\emptyset$.  
If $\pi:(X,D)\to (Y, \emptyset)$  is  a functorial   log  resolution of $Y$, then 
$$
\displaystyle\int_{X}c_{n}(T_{X}(-\log\, D)) =    \sum_{x\in  \Sing(v) } \mu_x(v),
  $$
  where $ \mu_x(v)$ is the Milnor number of  $v$ at $x$.
\end{cor}

\subsection*{Acknowledgments}
 The first named author   is partially supported by the Universit\`a degli Studi di Bari and by the
 PRIN 2022MWPMAB- "Interactions between Geometric Structures and Function Theories" and he is a member of INdAM-GNSAGA;
he was  partially supported by CNPq grant numbers 202374/2018-1, 400821/2016-8 and  Fapemig grant numbers  APQ-02674-21, APQ-00798-18,  APQ-00056-20; he is grateful to the University of Oxford for its hospitality during his visit, when the work on this project had begun.  
 The authors also thank anonymous referees for giving many suggestions that helped to improve the
presentation of the article.

\section{Preliminaries}

\subsection{Logarithmic  forms and vector fields}

\hyphenation{theo-ry}

Let $X$ be a complex manifold of dimension $n$ and let $D$ be a reduced hypersurface on $X$. Given a meromorphic $q$-form $\omega$ on $X$, we say that $\omega$ is a {\it logarithmic} $q$-form along $D$ at $x\in X$ if the following conditions occurs:

\begin{enumerate}
\item[$(i)$] $\omega$ is holomorphic on $X - D$;
\medskip
\item[$(ii)$] If $f=0$ is a reduced equation of $D$, locally at $x$, then $f \,\omega$ and $f\, d \omega$ are holomorphic.
\end{enumerate}
\noindent Denoting by $\Omega_{X,x}^q(\log\, D)$ the set of germs of logarithmic $q$-form  along $D$ at $x$, we define the following coherent sheaf of $\OO_X$-modules
$$
\Omega_{X}^q(\log\, D):= \bigcup_{x\in X} \Omega_{X,x}^q(\log\, D),
$$
which is called  by {\it sheaf of logarithmic $q$-forms along $D$}.  See \cite{Deligne}, \cite{Katz} and  \cite{Sai}  for details.

Now, given $x\in X$, let $v\in T_{X,x}$ be germ at $x$ of a holomorphic vector field on $X$. We say that $v$ is a {\it logarithmic vector field along of $D$ at $x$}, if $v$ satisfies the following condition: if $f=0$ is a equation of $D$, locally at $x$, then the derivation $v(f)$ belongs to the ideal $\langle f_x\rangle\OO_{X,x}$. Denoting by $T_{X,x}(-\log\, D)$ the set of germs of logarithmic vector field along of $D$ at $x$, we define the following coherent sheaf of $\OO_X$-modules
$$
T_{X}(-\log\, D):= \bigcup_{x\in X}T_{X,x}(-\log\, D),
$$
which is called  {\it sheaf of logarithmic vector fields along $D$}.  
 It follows from  \cite{Sai} that  if  $D$ is an analytic hypersurface with normal crossings  singularities, then  the sheaves $\Omega_{X}^1(\log\, D)$ and $T_{X}(-\log \,D)$ are  locally free, furthermore, the Poincar\'e residue map
$$
\mbox{Res}: \Omega_X^1(\log\, D) \longrightarrow \OO_{D} \cong \bigoplus_{i=1}^k\OO_{D_i}
$$
give the following exact sequence of sheaves on $X$
\begin{eqnarray}\label{exa_001}
0 \longrightarrow \Omega_X^1 \longrightarrow \Omega_X^1(\log\, D) \stackrel{ Res}{\longrightarrow} \bigoplus_{i=1}^k\OO_{D_i} \longrightarrow 0,
\end{eqnarray}
\noindent where $\Omega_X^1$ is the sheaf of holomorphic  $1$-forms on $X$ and $D_1,\ldots,D_k$ are the irreducible components of $D$.

\subsection{Singular one-dimensional holomorphic foliations }

A {\it one-dimensional holomorphic foliation $\F$} on $X$ is
given by the following data:
\begin{itemize}
  \item[$(i)$] an open covering $\mathcal{U}=\{U_{\alpha}\}$ of $X$;
  \item [$(ii)$] for each $U_{\alpha}$ a holomorphic vector field $v_\alpha \in TX|_{U_{\alpha}}$ ;
  \item [$(iii)$]for every non-empty intersection, $U_{\alpha}\cap U_{\beta} \neq \emptyset $, a
        holomorphic function $$g_{\alpha\beta} \in \mathcal{O}_X^*(U_\alpha\cap U_\beta);$$
\end{itemize}
such that $v_\alpha = g_{\alpha\beta}v_\beta$ in $U_\alpha\cap U_\beta$ and $g_{\alpha\beta}g_{\beta\gamma} = g_{\alpha\gamma}$ in $U_\alpha\cap U_\beta\cap U_\gamma$.

We denote by $K_{\F}$ the line bundle defined by the cocycle $\{g_{\alpha\beta}\}\in \mathrm{H}^1(X, \mathcal{O}^*)$. Thus, a one-dimensional holomorphic  foliation $\F$ on $X$ induces  a global holomorphic section $\vartheta_{\F}\in \mathrm{H}^0(X,T_X\otimes K_{\F})$.  
The line bundle $T_{\F}:= (K_{\F})^*  $ is called the  \emph{tangent  bundle} of $\F$. The singular set of $\F$ 
is $\singf=\{\vartheta_{\F}=0\}$. 

Throughout this paper  we will   assume that $\cod(\singf)\geq 2$.

Given $\F$ a one-dimensional holomorphic foliation on $X$, we say that $\F$ is \textit{logarithmic along $D$} if it satisfies the following condition: for all $x\in D - \Sing(D)$,  the vector $v_{\alpha}(x)$ belongs to $T_xD$, with $x\in U_{\alpha}$.


\subsubsection{Projective holomorphic  foliations  }
A foliation on a  complex projective space $\mathbb{P}^n$  is  called a \textit{projective foliation}. 
Let $\F $ be  a   projective foliation of dimension one   with tangent bundle $T_{\F}=\mathcal{O}_{\mathbb{P}^n}(r)$. The integer   $d:=r+1$ is called  the \textit{degree} of $\F$.
A projective foliation of dimension one of degree $d$ can be induced  by  a polynomial  vector field on $\C^{n+1}$ with homogeneous coefficients of degree $d$, see for instance \cite{CMS1,CMS2,Correa-toricas}.

\subsection{ Logarithmic and  homological indices  }
Let $D\subset X$  be a reduced hypersurface with  
a local equation $f\in \mathcal{O}_{X,x}$
in a neighborhood of a point $x\in D$.
Consider the  $\mathcal{O}_{D,x}$-module of
germs of regular  differentials of order $i$ on $D$ :
$$
\Omega_{D,x}^i=\frac{\Omega_{X,x}^i}{f\Omega_{X,x}^i+
df\wedge\Omega_{X,x}^{i-1}}.
$$
Now,  let  $\fol $ be a one-dimensional holomorphic foliation on $X$, with isolated
singularities, logarithmic along $D$. Let $x\in\singf$ be   and consider   a germ of vector field $v\in T_{X}(-\log \,D)|_{U}$ on $(U,x)$ tangent to $\fol $, where $U$ is  a neighborhood of  $x$.   Since $v$ is also  tangent to $(D,x)$ the interior multiplication $i_{v}$ induces the complex
$$
0\longrightarrow \Omega_{D,x}^{n-1} \stackrel{i_v}{\longrightarrow} \Omega_{D,x}^{n-2} \stackrel{i_v}{\longrightarrow}\cdots \stackrel{i_v}{\longrightarrow}\Omega_{D,x}^1 \stackrel{i_v}{\longrightarrow}\OO_{D,x}\longrightarrow0.
$$ 
The \emph{homological index} is defined as the Euler characteristic of the complex $(\Omega_{D,x}^{\bullet}, i_v)$
$$
Ind_{hom}(\fol,D,x)=\sum_{i=0}^{n-1}(-1)^i\dim H_i(\Omega_{D,x}^{\bullet}, i_v). 
$$
Since  the vector field $v$ has an isolated singularity at  $x$, then the $i_v$-homology groups
of the complex $\Omega_{D,x}^{\bullet}$ are finite-dimensional vector spaces and   the Euler characteristic is well defined. Similarly, the Euler characteristic of the Koszul complex 
  $(\Omega_{X,x}^{\bullet}, i_v)$ associated to $v$  is well defined and the \emph{ Milnor number}  of $\fol$ at $x$ is defined  by 
$$
\mu(\fol,x)=\sum_{i=0}^{n}(-1)^i\dim H_i(\Omega_{X,x}^{\bullet}, i_v). 
$$
We recall that  Milnor number of $\fol$ at $x$   coincides with the  Poincar\'e Hopf   index of $v$ at $x$. For more details on  residues of holomorphic foliations  see \cite{BB1, BarSaeSuw,  C-L,CPS, Suw2}. 

The homological index was  introduced by G\'omez-Mont in \cite{a08} and it coincides  with the  GSV-index  introduced by  G\'omez-Mont,   Seade and  Verjovsky in  \cite{GSV}.    The concept of GSV-index has been  extended to  more general contexts, we refer to the works \cite{SeaSuw1,a07, a03, BarSaeSuw, Suw2, CoMa}.

Aleksandrov  has proved in \cite[Proposition 1]{Ale_p} that the  following formula  holds 
\begin{eqnarray}\label{eq001}
Log(\fol,D,x) = \mu_x(\fol) - Ind_{hom}(\fol,D,x),
\end{eqnarray}
where $\mu_x(\fol)$ is  the Milnor number    of $\fol$ at $x$.
Denote by  $D_{reg}:= D \setminus \Sing(D)$  the regular part of $D$.
If $x\in \Sing(\fol)\cap D_{reg}$    then
\begin{eqnarray}\label{eq0_02}
Log(\fol,D,x) = \mu_x(\fol) - \mu_x(\fol|_{{D}_{reg}}),
\end{eqnarray}
\noindent  since in this case the GSV-index of $\fol$ at $x$ coincides with the Milnor number of  the restriction $\fol|_{{D}_{reg}}$ at $x$. Furthermore, 
\begin{eqnarray}\nonumber
Log(\fol,D,x)= 0,
\end{eqnarray}
\noindent whenever $x$ is a non-degenerate singularity of $\fol$.

\hyphenation{using}

\hyphenation{sin-gu-la-ri-ties}
\hyphenation{sin-gu-la-ri-ty}

\hyphenation{fo-lia-tion}

\section{Proof of Theorem \ref{teo1_1}} \label{sec03}
 If $x\in \Sing(\fol)\cap D_{reg}$,  then it follows from (\ref{eq0_02}) that 
$$
Log(\fol,D,x) = \mu_x(\fol) - \mu_x(\fol|_{{D}_{reg}}).
$$
 In order to simplify the proof we also will adopt the notation
 \begin{eqnarray} \label{eq0_03}
Log(\fol,D,x) = \mu_x(\fol),
\end{eqnarray}
\noindent  if    $x\in X\setminus D$.

Let $D = D_ {1} \cup \ldots \cup D_{k}$ be the decomposition of $D$ into irreducible components. Fixing an irreducible component, let us say $D_k$, we define $ \widehat{D}_k:= D_1 \cup \ldots \cup D_{k-1}$. The intersection $\widehat{D}_k \cap D_k = (D_1 \cap D_k)\, \cup \ldots \cup\, (D_{k-1} \cap D_k)$ is an analytic hypersurface with normal crossings  singularities on the $(n-1)$-dimensional smooth submanifold  $D_k$, having $k-1$ irreducible components. Furthermore, if $\fol $ is logarithmic along $D$, then its restriction $\fol|_{D_k}$ on $D_k$  is logarithmic along $\widehat{D}_k \cap D_k$. 

\begin{lema}\label{lemanovo0} Let $X$, $D = D_ {1} \cup \ldots \cup D_{k}$ and $\fol$ be as described above. Then, for all $x\in \widehat{D}_k \cap D_k$, we have 
\begin{eqnarray} \label{eq0_04}
Log (\fol,D,x )= \displaystyle Log (\fol,\widehat{D}_k,x) - Log (\fol|_{D_k},\widehat{D}_k \cap D_k,x ).   
\end{eqnarray}
\end{lema}
\begin{demostracion}
Given $x\in \widehat{D}_k \cap D_k$,   let $(z_1,\ldots,z_n)$ be 
a local coordinate system in a neighborhood of $x$,  such that each $D_j$ is locally defined by equation  $z_j = 0$, so $D$ is defined by the equation $z_1\cdot \cdot  \cdot z_k = 0$.
 Moreover, we have that the set
$$
\displaystyle \left\{\frac{dz_1}{z_1},\ldots, \frac{dz_k}{z_k},dz_{k+1},\ldots,dz_{n}\right\}
$$
\noindent constitutes a system of $\OO_{X,x}$ - free basis for $\Omega_{X,x}^1(\log\, D)$. In general, for all $q = 1,2,\ldots, n$, we have that 
$$
\Omega_{X,x}^q(\log\, D) = \sum_{i_1,\ldots,i_q}\OO_{X,x}\, \omega_{i_1}\wedge\ldots\wedge \omega_{i_q}, 
$$ 
\noindent where $\omega_1 = \displaystyle \frac{dz_1}{z_1}\,\,,\ldots,\,\, \omega_k = \frac{dz_k}{z_k},
\,\, \omega_{k+1} = dz_{k+1} \,\,,\ldots,\,\, \omega_{n} = dz_{n}$, see  \cite[pg. 270]{Sai}.  Thus, we have the map of complexes 
$$
\Omega_{X,x}^{\bullet}(\log\,D) \displaystyle\stackrel{\cdot z_k}{\longrightarrow} \Omega_{X,x}^{\bullet}(\log\, \widehat{D}_k)
$$
\noindent which defines the following exact sequence
\begin{eqnarray}\nonumber
\displaystyle 0 \longrightarrow \left(\Omega_{X,x}^{\bullet}(\log\,D)\right) \displaystyle \stackrel{\cdot z_k}{\longrightarrow} \left(\Omega_{X,x}^{\bullet}(\log\, \widehat{D}_k)\right) \stackrel{i^{\ast}}{\longrightarrow} \left(\Omega_{D_1,x}^{\bullet}(\log\, (\widehat{D}_k\cap D_k))\right) \longrightarrow 0,
\end{eqnarray} 
\noindent where $i^{\ast}$ denotes the pullback of inclusion map $i: D_k \hookrightarrow X $.
Therefore,  the result follows from the  additivity of the holomorphic Euler characteristic. \\ 
$\square$
\end{demostracion}

\begin{lema} \label{lemanovo} Let $X$, $D = D_ {1} \cup \ldots \cup D_{k}$ and $\fol$ be as described above. Then, for all $x\in  \widehat{D}_k -  (\widehat{D}_k \cap D_k)$, we have that 
\begin{eqnarray} \label{eq0_04}
Ind_{hom} (\fol,\widehat{D}_k ,x )= Ind_{hom} (\fol,D,x ).
\end{eqnarray}
\end{lema}
\begin{demostracion}
Given $x\in D$, since $D$ is a hypersurface with normal crossings  singularities, we can choose  a local coordinate system in a neighborhood of  $x$ such that each $D_j$ is locally defined by the  equation $z_j = 0$. If $x\in  \widehat{D}_k -  (\widehat{D}_k \cap D_k)$, then $\widehat{D}_k$ and $D$ are both defined by the same equation
$
z_1\cdots z_{k-1} = 0.
$
\noindent Thus, for all $q = 0,1,\ldots,n$, we have that 
$
\Omega^q_{\widehat{D}_k,x} \cong \Omega^q_{D,x},
$
 and, consequently, we obtain
 \begin{eqnarray} \nonumber
Ind_{hom} (\fol,\widehat{D}_k ,x )&=& \displaystyle\sum_{i=0}^{n}(-1)^i dim H_i(\Omega^\bullet_{\widehat{D}_k,x} ,i_{v|_{\widehat{D}_k}})  \\\nonumber \\\nonumber  &=& \displaystyle\sum_{i=0}^{n}(-1)^i dim H_i(\Omega^\bullet_{D,x} ,i_{v|_{D}})
\\\nonumber \\\nonumber &=&  Ind_{hom} (\fol,D,x ).
\end{eqnarray}
$\square$
\end{demostracion}

We will use the following multiple index notation: for each multi-index $J = (j_1,\ldots,j_k)$ and $J'=(j_1',\ldots,j_{k-1}')$, with $1\leq j_l,j'_t\leq n$, we denote 
\begin{eqnarray}\nonumber
c_1(D)^J & = & c_1([D_1])^{j_1}\cdots c_1([D_k])^{j_k}, \\\nonumber
c_1(\widehat{D}_k)^{J'} &=& c_1([D_1])^{j'_1}\cdots c_1([D_{k-1}])^{j'_{k-1}}.
\end{eqnarray}

\begin{lema} \label{p13} In the above conditions, for each $i = 1,\ldots, n$, we have
\begin{eqnarray}\nonumber
\displaystyle c_{i}(\Omega_{X}^1(\log\, D)) = \sum^{i}_{l=0}\sum_{\mid J \mid = l} c_{i-l}(\Omega_X^1)c_1(D)^{J}.
\end{eqnarray}
\end{lema}
\begin{demostracion} By using the exact sequence (\ref{exa_001}), we get 
\begin{eqnarray} \nonumber
c_i(\Omega_X^1(\log\, D)) &=& \sum^{i}_{l=0}c_{i-l}(\Omega^1_X)c_l\left(\bigoplus_{i=1}^k\OO_{D_i}\right) \\\nonumber \\\nonumber  &=& \sum^{i}_{l=0}c_{i-l}(\Omega^1_X)\left (\sum_{j_1+\ldots +j_k=l}c_{j_1}(\OO_{D_1})\ldots c_{j_k}(\OO_{D_k})\right) \\\nonumber \\\nonumber &=&\sum^{i}_{l=0}c_{i-l}(\Omega_X^1)\left(\sum_{j_1+\ldots+j_k = l}c_1([D_1])^{j_1}\ldots c_1([D_k])^{j_k}\right),
\end{eqnarray}
\noindent where in the  last equality we use    the   relations
$$
c_j(\OO_{D_i}) = c_j(\OO_X - \OO(-D_i)) =  c_1([D_i])^j,\,\,\,\, j = 1,\ldots,n.
$$
$\square$
\end{demostracion}

\begin{lema} \label{p1300} In the above conditions, for each irreducible  component $D_j$, we have that 
\begin{eqnarray}\nonumber
\displaystyle c_{i}(\Omega_{X}^1)|_{D_j} = c_{i}(\Omega_{D_j}^1) - c_{i-1}(\Omega_{D_j}^1)c_i([D_j])|_{D_j},\,\,\,\forall i = 1,\ldots, n-1.
\end{eqnarray}
\end{lema}
\begin{demostracion} 
This  follows  by taking the total  Chern class in the exact sequence
$$
0\to  T_{D_j}  \to T_X|_{D_j}  \to  [D_j]|_{D_j}\to 0.
$$
$\square$
\end{demostracion}

\begin{remark}\label{ref-virtual}
We recall  that the  total Chern class of the virtual vector bundle $T_X(-\log\, D) -  L$ is defined   by 
$$
c(T_X(-\log\, D) -  L)=\frac{c(T_X(-\log\, D))}{c(L)} 
$$
and the  $n$-th Chern class of $T_X(-\log\, D) -  L$ is by definition the component of $\frac{c(T_X(-\log\, D))}{c(L)} $ in dimension $2n$, see \cite[section 4]{BB1}. Moreover, we have that 
$$
c_n(T_X(-\log\, D) -  L)= c_n(T_X(-\log\, D) \otimes L^*). 
$$
\end{remark}

\begin{lema} \label{p0013} In the above conditions, if $L$ is a holomorphic line bundle  on $X$, then
the following   holds: 
\begin{eqnarray} \label{rel001}
\int_{X}c_n(T_X(-\log\, D) -  L) =\\\nonumber \sum_{j=0}^n\sum^{n-j}_{l=0} \sum_{\mid J \mid = l}  \int_{X} (-1)^{n-j}c_{n-j-l}(\Omega_X^1) c_1(D)^{J}c_1(L^{\ast})^j.
\end{eqnarray}
\noindent In particular, we have that 
\begin{eqnarray} \label{rel002}
\int_{X}c_n(T_X(-\log\, \widehat{D}_k) -  L) =\\\nonumber \sum_{j=0}^n\sum^{n-j}_{l=0} \sum_{\mid J' \mid = l}  \int_{X} (-1)^{n-j}c_{n-j-l}(\Omega_X^1) c_1(\widehat{D}_k)^{J'}c_1(L^{\ast})^j.
\end{eqnarray}
\noindent and
\begin{eqnarray} \label{rel003}
\int_{D_k}c_{n-1}(T_{D_k}(-\log\, (\widehat{D}_k\cap D_k)) - L|_{D_k}) = \\\nonumber \sum_{j=0}^{n-1}\sum^{n-1-j}_{l=0} \sum_{\mid J' \mid = l}  \int_{D_k} (-1)^{n-1-j}c_{n-1-j-l}(\Omega_{D_k}^1) c_1(\widehat{D}_k)^{J'}c_1(L^{\ast})^j.
\end{eqnarray} 
\end{lema}

\begin{demostracion} On the one hand, since by remark \ref{ref-virtual} we have that 
$$
c_n(T_X(-\log\, D) - L)= c_n(T_X(-\log\, D) \otimes L^*)= \sum_{j=0}^n  c_{n-j}(T_{X}(-\log\, D))c_1(L^{\ast})^j 
$$
and $ c_{n-j}(T_{X}(-\log\, D))= (-1)^{n-j} c_{n-j}(\Omega_{X}^1(\log\, D))$ 
 we get
\begin{eqnarray}\nonumber
\displaystyle \int_{X}c_n(T_X(-\log\, D) - L)  
&=& \int_{X}\sum_{j=0}^n(-1)^{n-j} c_{n-j}(\Omega_{X}^1(\log\, D))c_1(L^{\ast})^j.
\end{eqnarray}
\noindent On the other hand,  by Lemma \ref{p13}, we obtain 
$$
c_{n-j}(\Omega_{X}^1(\log \, D)) = \sum^{n-j}_{l=0} \sum_{\mid J \mid = l} c_{n-j-l}(\Omega_X^1)c_1(D)^{J}.
$$ 
 Substituting this, we obtain (\ref{rel001}). We get 
the relation (\ref{rel002})  by taking $D=\widehat{D}_k$ in the  relation (\ref{rel001}). Analogously,   applying  the relation (\ref{rel002}), we can obtain (\ref{rel003}) by taking $X = D_k$ as a complex manifold of dimension $n-1$ and $D = \widehat{D}_k \cap D_k$ as an analytic subvariety of  $D_k$ with normal crossings  singularities. $\square$
\end{demostracion}

\begin{prop}\label{lema4}
In the above conditions, if $L$  is a holomorphic line bundle  on $X$, then
{\small {\begin{eqnarray} \nonumber
\int_{X}c_n(T_X(-\log\, D) -  L) &=& \int_{X}c_n(T_X(-\log (\widehat{D}_k) - L) - \int_{D_k}c_{n-1}(T_{D_k} 
(- \log( \widehat{D}_k \cap D_k) - L|_{D_k}).
\end{eqnarray}}}

\noindent In particular, if $D$ is  irreducible  we get 
\begin{eqnarray}\label{eq050}
\displaystyle\int_{X}c_{n}(T_{X}(-\log\, D)- L) &=& \int_{X}c_{n}(T_{X}- L) - \int_{D}c_{n-1}(T_{D}- L|_{D}).
\end{eqnarray}
\end{prop}

\begin{demostracion} By Lemma \ref{p0013},   it is enough to show that the following equality occurs

\begin{eqnarray}\nonumber
\sum_{j=0}^n\sum^{n-j}_{l=0} \sum_{\mid J \mid = l}  \int_{X} (-1)^{n-j}c_{n-j-l}(\Omega_X^1) c_1(D)^{J}c_1(L^{\ast})^j &=&\\\nonumber 
\sum_{j=0}^n\sum^{n-j}_{l=0} \sum_{\mid J' \mid = l}  \int_{X} (-1)^{n-j}c_{n-j-l}(\Omega_X^1) c_1(\widehat{D}_k)^{J'}c_1(L^{\ast})^j &-&\\\nonumber - \sum_{j=0}^{n-1}\sum^{n-1-j}_{l=0} \sum_{\mid J' \mid = l}  \int_{D_k} (-1)^{n-1-j}c_{n-1-j-l}(\Omega_{D_k}^1) c_1(\widehat{D}_k)^{J'}c_1(L^{\ast})^j.
\end{eqnarray}

\noindent Indeed, we can decompose the sum on the left hand side into the terms with $l=0$ and those with $l\geq 1$ as follows:

\begin{eqnarray}\label{eq00003}
\sum_{j=0}^n\sum^{n-j}_{l=0} \sum_{\mid J \mid = l}  \int_{X} (-1)^{n-j}c_{n-j-l}(\Omega_X^1) c_1(D)^{J}c_1(L^{\ast})^j &=&  \\\nonumber\\\nonumber\\\nonumber  =  \displaystyle \sum^{n}_{j=0}\int_{X}(-1)^{n-j}c_{n-j}(\Omega^1_X)c_1(L^{\ast})^j + \sum_{j=0}^{n-1}\sum^{n-j}_{l=1} \sum_{\mid J \mid = l}  \int_{X} (-1)^{n-j}c_{n-j-l}(\Omega_X^1) c_1(D)^{J}c_1(L^{\ast})^j.
\end{eqnarray}
\noindent The second sum on the right hand side can   be computed as follows:  
\begin{eqnarray}\nonumber
\sum_{j=0}^{n-1}\sum^{n-j}_{l=1} \sum_{\mid J \mid = l}  \int_{X} (-1)^{n-j}c_{n-j-l}(\Omega_X^1) c_1(D)^{J}c_1(L^{\ast})^j = \\\nonumber\\\nonumber
\sum_{j=0}^{n-1}\sum^{n-j}_{l=1} \sum_{\mid J' \mid = l}  \int_{X} (-1)^{n-j}c_{n-j-l}(\Omega_X^1) c_1(\widehat{D}_k)^{J'}c_1(L^{\ast})^j + \\\nonumber\\\nonumber  + \sum_{j=0}^{n-1}\sum^{n-j}_{l=1} \sum_{\substack{\mid J \mid = l \\ j_k \geq 1}}  \int_{X} (-1)^{n-j} c_{n-j-l}(\Omega_X^1) c_1([D_1])^{j_1}\ldots c_1([D_k])^{j_k}c_1(L^{\ast})^j.
\end{eqnarray}
\noindent  By using  that $c_1 ([D_k])$ is the  Poincar\'e dual to the fundamental class of $D_k$,
we obtain:
\begin{eqnarray}\nonumber
\sum_{j=0}^{n-1}\sum^{n-j}_{l=1} \sum_{\mid J \mid = l}  \int_{X} (-1)^{n-j}c_{n-j-l}(\Omega_X^1) c_1(D)^{J}c_1(L^{\ast})^j = \\\nonumber\\\nonumber
\sum_{j=0}^{n-1}\sum^{n-j}_{l=1} \sum_{\mid J' \mid = l}  \int_{X} (-1)^{n-j}c_{n-j-l}(\Omega_X^1) c_1(\widehat{D}_k)^{J'}c_1(L^{\ast})^j + \\\nonumber\\\nonumber + \sum_{j=0}^{n-1}\sum^{n-j}_{l=1} \sum_{\substack{\mid J \mid = l \\ j_k \geq 1}}  \int_{D_k} (-1)^{n-j} c_{n-j-l}(\Omega_X^1) c_1([D_1])^{j_1}\ldots c_1([D_k])^{j_k - 1}c_1(L^{\ast})^j.
\end{eqnarray}

\noindent Now, by using the relation of Lemma \ref{p1300} we get

\begin{eqnarray}\nonumber
\sum_{j=0}^{n-1}\sum^{n-j}_{l=1} \sum_{\substack{\mid J \mid = l \\ j_k \geq 1}}  \int_{D_k} (-1)^{n-j} c_{n-j-l}(\Omega_X^1) c_1([D_1])^{j_1}\ldots c_1([D_k])^{j_k - 1}c_1(L^{\ast})^j = \\\nonumber\\\nonumber = \sum_{j=0}^{n-1}\sum^{n-j}_{l=1} \sum_{\substack{\mid J \mid = l \\ j_k \geq 1}}  \int_{D_k} (-1)^{n-j} c_{n-j-l}(\Omega_{D_k}^1) c_1([D_1])^{j_1}\ldots c_1([D_k])^{j_k - 1}c_1(L^{\ast})^j - \\\nonumber\\\nonumber  - \sum_{j=0}^{n-1}\sum^{n-j-1}_{l=1} \sum_{\substack{\mid J \mid = l \\ j_k \geq 1}}  \int_{D_k} (-1)^{n-j} c_{n-j-1-l}(\Omega_{D_k}^1) 
c_1(D)^{J}c_1(L^{\ast})^j = \\\nonumber\\\nonumber  = \sum_{j=0}^{n-1}  \int_{D_k} (-1)^{n-j} c_{n-j-1}(\Omega_{D_k}^1) c_1(L^{\ast})^j + \sum_{j=0}^{n-1}\sum^{n-j-1}_{l=1} \sum_{\mid J' \mid = l}  \int_{D_k} (-1)^{n-j} c_{n-j-1-l}(\Omega_{D_k}^1) c_1(\widehat {D}_k)^{J'}c_1(L^{\ast})^j = \\\nonumber\\\nonumber = - \sum_{j=0}^{n-1}\sum^{n-1-j}_{l=0} \sum_{\mid J' \mid = l}  \int_{D_k} (-1)^{n-1-j} c_{n-1-j-l}(\Omega_{D_k}^1) c_1(\widehat {D}_k)^{J'}c_1(L^{\ast})^j. 
\end{eqnarray}
Hence,
\begin{eqnarray}\nonumber
\sum_{j=0}^{n-1}\sum^{n-j}_{l=1} \sum_{\mid J \mid = l}  \int_{X} (-1)^{n-j}c_{n-j-l}(\Omega_X^1) c_1(D)^{J}c_1(L^{\ast})^j = \\\nonumber\\\nonumber
\sum_{j=0}^{n-1}\sum^{n-j}_{l=1} \sum_{\mid J' \mid = l}  \int_{X} (-1)^{n-j}c_{n-j-l}(\Omega_X^1) c_1(\widehat{D}_k)^{J'}c_1(L^{\ast})^j - \\\nonumber\\\nonumber - \sum_{j=0}^{n-1}\sum^{n-1-j}_{l=0} \sum_{\mid J' \mid = l}  \int_{D_k} (-1)^{n-1-j} c_{n-1-j-l}(\Omega_{D_k}^1) c_1(\widehat {D}_k)^{J'}c_1(L^{\ast})^j,
\end{eqnarray}

\noindent and we complete  the calculation of the second sum.
Replacing it in the initial equality (\ref{eq00003}), we obtain

\begin{eqnarray}\nonumber 
\sum_{j=0}^n\sum^{n-j}_{l=0} \sum_{\mid J \mid = l}  \int_{X} (-1)^{n-j}c_{n-j-l}(\Omega_X^1) c_1(D)^{J}c_1(L^{\ast})^j &=&  \\\nonumber\\\nonumber  =  \displaystyle \sum^{n}_{j=0}\int_{X}(-1)^{n-j}c_{n-j}(\Omega^1_X)c_1(L^{\ast})^j + \sum_{j=0}^{n-1}\sum^{n-j}_{l=1} \sum_{\mid J' \mid = l}  \int_{X} (-1)^{n-j}c_{n-j-l}(\Omega_X^1) c_1(\widehat{D}_k)^{J'}c_1(L^{\ast})^j - \\\nonumber\\\nonumber - \sum_{j=0}^{n-1}\sum^{n-1-j}_{l=0} \sum_{\mid J' \mid = l}  \int_{D_k} (-1)^{n-1-j} c_{n-1-j-l}(\Omega_{D_k}^1) c_1(\widehat{D}_k)^{J'}c_1(L^{\ast})^j = \\\nonumber\\\nonumber  = \sum_{j=0}^n\sum^{n-j}_{l=0} \sum_{\mid J' \mid = l}  \int_{X} (-1)^{n-j}c_{n-j-l}(\Omega_X^1) c_1(\widehat{D}_k)^{J'}c_1(L^{\ast})^j - \\\nonumber\\\nonumber - \sum_{j=0}^{n-1}\sum^{n-1-j}_{l=0} \sum_{\mid J' \mid = l}  \int_{D_k} (-1)^{n-1-j}c_{n-1-j-l}(\Omega_{D_k}^1) c_1(\widehat{D}_k)^{J'}c_1(L^{\ast})^j .
\end{eqnarray}

Now, if $D$ is  irreducible, we can repeat the same argument that we have used above  in order  to obtain the following formula: 
\begin{eqnarray} \nonumber
\sum_{j=0}^n\sum^{n-j}_{l=0}   \int_{X} (-1)^{n-j}c_{n-j-l}(\Omega_X^1) c_1([D])^{l}c_1(L^{\ast})^j = \sum_{j=0}^n  \int_{X} (-1)^{n-j}c_{n-j}(\Omega_X^1) c_1(L^{\ast})^j - \\\nonumber \\\nonumber \displaystyle -\sum_{j=0}^{n-1}\int_{D} (-1)^{n-1-j}c_{n-1-j}(\Omega_{D}^1) c_1(L^{\ast})^j.
\end{eqnarray} 

\noindent Thus, we get

\begin{eqnarray} \nonumber
\int_{X}c_n(T_X(-\log\, D) -  L) &=& \sum_{j=0}^n  \int_{X} (-1)^{n-j}c_{n-j}(\Omega_X^1) c_1(L^{\ast})^j - \sum_{j=0}^{n-1}\int_{D} (-1)^{n-1-j}c_{n-1-j}(\Omega_{D}^1) c_1(L^{\ast})^j\\\nonumber \\\nonumber &=& \sum_{j=0}^n  \int_{X} c_{n-j}(T_X) c_1(L^{\ast})^j - \sum_{j=0}^{n-1}\int_{D} c_{n-1-j}(T_{D}) c_1(L^{\ast})^j\\\nonumber \\\nonumber &=& \int_{X}c_{n}(T_{X}- L) - \int_{D}c_{n-1}(T_{D}- L|_{D}).
\end{eqnarray} 
$\square$
\end{demostracion}

In order to prove  the Theorem \ref{teo1_1}  we will use  the  induction principle  on the number of irreducible components of $D$. Indeed, if the number of irreducible component of $D$ is 1, then we can invoke the formula (\ref{eq050}) of the Proposition \ref{lema4} to obtain the following 
\begin{eqnarray}\nonumber
\displaystyle\int_{X}c_{n}(T_{X}(-\log\, D)- T_{\fol}) &=& \int_{X}c_{n}(T_{X}- T_{\fol}) - \int_{D}c_{n-1}(T_{D}- T_{\fol})\\\nonumber\\\nonumber
&=& \sum_{x\in \Sing(\fol) \cap X} \mu_x(\fol) - \sum_{x\in \Sing(\fol) \cap D} \mu_x(\fol|_{D}), 
\end{eqnarray}
\noindent where in the last step  we apply the Baum-Bott classical formula, see \cite{BB1}. 
By using the disjoint decomposition $X = (X - D) \cup D$, we get
\begin{eqnarray}\nonumber
\displaystyle\int_{X}c_{n}(T_{X}(-\log\, D)- T_{\fol})&=& \sum_{x\in   \Sing(\fol) \cap (X - D) } \mu_x(\fol) + \\\nonumber \\\nonumber &+& \left ( \sum_{x\in  \Sing(\fol) \cap (D) } \mu_x(\fol) - \sum_{x\in  \Sing(\fol) \cap D  } \mu_x(\fol|_{D})\right)
\end{eqnarray}

\noindent But, by relation (\ref{eq0_03}), we have

\begin{eqnarray}\nonumber
\sum_{x\in  \Sing(\fol) \cap (X - D) } Log(\fol,D, x) = \sum_{x\in   \Sing(\fol) \cap (X - D) } \mu_x(\fol).
\end{eqnarray}

\noindent  On the other hand, since $D$ smooth, it follows from equality (\ref{eq0_02}) that

\begin{eqnarray}\nonumber
\sum_{x\in  \Sing(\fol) \cap (D) } Log(\fol,D, x) = \sum_{x\in  \Sing(\fol) \cap (D)  } \mu_x(\fol) - \sum_{x\in \Sing(\fol) \cap {D}  } \mu_x(\fol|_{D}).
\end{eqnarray}

\noindent Therefore, we get 

\begin{eqnarray}\nonumber
\displaystyle\int_{X}c_{n}(T_{X}(-\log\, D)- T_{\fol}) = \sum_{x\in Sing(\fol)} Log(\fol,D, x).
\end{eqnarray}

\hyphenation{hy-po-the-sis}

Let us suppose that for every analytic hypersurface on $X$, satisfying the hypothesis of Theorem \ref{teo1_1}  and having $k-1$ irreducible components, the formula of Theorem \ref{teo1_1} holds. Let $D$ be an analytic hypersurface on $X$ with $k$ irreducible components, satisfying the hypotheses of the Theorem \ref{teo1_1}. We will prove that the formula   of Theorem \ref{teo1_1} holds for $D$.

We know that $\widehat{D}_k$ is an analytic hypersurface on $X$ and  $\widehat{D}_k \cap {D}_k$ is an analytic hypersurface on $D_k$, both with normal crossings  singularities and having exactly $k - 1$ irreducible components. Moreover, $\fol$ and its restriction  $\fol|_{D_k}$ on $D_k$ are logarithmic along $D_k$ and $\widehat{D}_k \cap {D}_k$, respectively. Thus, we can use the induction hypothesis and we obtain
\begin{eqnarray}\nonumber 
 \int_{X}c_n(T_X(-\log\, \widehat{D}_k)- T_{\fol}) = \sum_{x\in \Sing(\fol)} Log(\fol,\widehat{D}_k, x)
\end{eqnarray}

\noindent and

\begin{eqnarray}\nonumber
\int_{D_{k}}c_{n-1}(T_{D_k}(-\log\,(\widehat{D}_k \cap {D}_k))- T_{\fol}|_{D_{k}}) = \sum_{x\in \Sing(\fol|_{D_k})}  Log(\fol|_{D_k},\widehat{D}_k \cap {D}_k, x).
\end{eqnarray}

By  the Proposition \ref{lema4}, we get
\begin{eqnarray}\nonumber
\displaystyle\int_{X}c_{n}(T_{X}(-\log\, D)- T_{\fol}) = \sum_{x\in \Sing(\fol)} Log(\fol,\widehat{D}_k, x) - \sum_{x\in \Sing(\fol|_{D_k})} Log(\fol|_{D_k},\widehat{D}_k \cap {D}_k, x).
\end{eqnarray}

\noindent Thus,  it is enough to show that the following equality holds 
$$
\sum_{x\in \Sing(\fol)} Log(\fol,\widehat{D}_k, x) = \sum_{x\in \Sing(\fol)} Log(\fol,D, x) + \sum_{x\in \Sing(\fol|_{D_k})}  Log(\fol|_{D_k},\widehat{D}_k \cap {D}_k, x). 
$$

\noindent Indeed, by  using the following disjoint decomposition
$$
X = (X - D) \cup \left\{ \widehat{D}_k \cap \left[ D - \left(\widehat{D}_k \cap D_k\right)\right]\right\} \cup \left[ \left(D - \widehat{D}_k \right) \cap D_k \right] \cup \left( \widehat{D}_k \cap D_k\right),
$$

\noindent we get

\begin{eqnarray}\nonumber
\sum_{x\in \Sing(\fol)} Log(\fol,\widehat{D}_k, x) &=& \sum_{\substack{x\in \Sing(\fol),\\x\in (X - D)}} Log(\fol,\widehat{D}_k, x) + \sum_{\substack{x\in \Sing(\fol),\\x\in \widehat{D}_k \cap \left[ D - \left(\widehat{D}_k \cap D_k\right)\right]}} Log(\fol,\widehat{D}_k, x) +\\\nonumber\\\nonumber &+& \sum_{\substack{x\in \Sing(\fol),\\x\in \left(D - \widehat{D}_k \right) \cap D_k }}
 Log(\fol,\widehat{D}_k, x) + \sum_{\substack{x\in \Sing(\fol),\\x\in \widehat{D}_k \cap D_k}}
 Log(\fol,\widehat{D}_k, x). 
\end{eqnarray}

\noindent We can rewrite each of the above sums in a more appropriate way, as follows: in the first one by using the relation (\ref{eq0_03})  we obtain
$$
\sum_{\substack{x\in \Sing(\fol),\\x\in ( X - D) }} Log(\fol,\widehat{D}_k, x) = \sum_{\substack{x\in \Sing(\fol),\\x\in( X - D)}} \mu_x(\fol) = \sum_{\substack{x\in \Sing(\fol),\\x\in (X - D)}} Log(\fol,D, x).
$$
\noindent In the second sum by  using the relation (\ref{eq001}) and the Lemma \ref{lemanovo}  we get
\begin{eqnarray}\nonumber
\sum_{\substack{x\in \Sing(\fol),\\x\in \widehat{D}_k \cap \left[ D - \left(\widehat{D}_k \cap D_k\right)\right]}} Log(\fol,\widehat{D}_k, x) &=&  \sum_{\substack{x\in \Sing(\fol),\\x\in \widehat{D}_k \cap \left[ D - \left(\widehat{D}_k \cap D_k\right)\right]}} \left[ \mu_x(\fol) - Ind_{hom} (\fol,\widehat{D}_k, x)\right]  \\\nonumber\\\nonumber &=& \sum_{\substack{x\in \Sing(\fol),\\x\in \widehat{D}_k \cap \left[ D - \left(\widehat{D}_k \cap D_k\right)\right]}} \left[ \mu_x(\fol) - Ind_{hom} (\fol,D, x)\right]  \\\nonumber\\\nonumber &=& \sum_{\substack{x\in \Sing(\fol),\\x\in \widehat{D}_k \cap \left[ D - \left(\widehat{D}_k \cap D_k\right)\right]}} Log(\fol,D, x).  
\end{eqnarray}
\noindent In the next sum by  using the relation (\ref{eq0_03}) again  we have 
\begin{eqnarray}\nonumber
\sum_{\substack{x\in \Sing(\fol),\\x\in \left(D - \widehat{D}_k \right) \cap D_k }}
 Log(\fol,\widehat{D}_k, x) &=&  \sum_{\substack{x\in \Sing(\fol),\\x\in \left(D - \widehat{D}_k \right) \cap D_k }}
 \mu_x(\fol)  \\\nonumber\\\nonumber &=&  \sum_{\substack{x\in \Sing(\fol),\\x\in \left(D - \widehat{D}_k \right) \cap D_k }}
 \mu_x(\fol) - \sum_{\substack{x\in \Sing(\fol),\\x\in \left(D - \widehat{D}_k \right) \cap D_k }}
 \mu_x(\fol|_{D}) + \sum_{\substack{x\in \Sing(\fol),\\x\in \left(D - \widehat{D}_k \right) \cap D_k }}
 \mu_x(\fol|_{D}), 
\end{eqnarray}
\noindent and, since $\left(D - \widehat{D}_k \right) \cap D_k \subset D_{reg}$, the identity  (\ref{eq0_02}) gives  us
\begin{eqnarray}\nonumber
\sum_{\substack{x\in \Sing(\fol),\\x\in \left(D - \widehat{D}_k \right) \cap D_k }}
 Log(\fol,\widehat{D}_k, x) &=&  \sum_{\substack{x\in \Sing(\fol),\\x\in \left(D - \widehat{D}_k \right) \cap D_k }}
 Log(\fol,D, x) + \sum_{\substack{x\in \Sing(\fol),\\x\in \left(D - \widehat{D}_k \right) \cap D_k }}
 \mu_x(\fol|_{D}). 
\end{eqnarray}
\noindent Since $\fol|_{D}$ and $\fol|_{D_k}$ coincide locally around   $x\in \left(D - \widehat{D}_k \right) \cap D_k$, we obtain
\begin{eqnarray}\nonumber
\sum_{\substack{x\in \Sing(\fol),\\x\in \left(D - \widehat{D}_k \right) \cap D_k }}
 Log(\fol,\widehat{D}_k, x) = \sum_{\substack{x\in \Sing(\fol),\\x\in \left(D - \widehat{D}_k \right) \cap D_k }}
 Log(\fol,D, x) + \sum_{\substack{x\in \Sing(\fol),\\x\in \left(D - \widehat{D}_k \right) \cap D_k }}
 \mu_x(\fol|_{D_k}).
\end{eqnarray}
\noindent Finally, in the last sum,  by    Lemma \ref{lemanovo0} we get that 
$$
\sum_{\substack{x\in \Sing(\fol),\\x\in \widehat{D}_k \cap D_k}}
 Log(\fol,\widehat{D}_k, x) = \sum_{\substack{x\in \Sing(\fol),\\x\in \widehat{D}_k \cap D_k}}
 Log(\fol,D, x) + \sum_{\substack{x\in \Sing(\fol),\\x\in \widehat{D}_k \cap D_k}}
 Log(\fol|_{D_k},\widehat{D}_k \cap D_k, x).
$$
Now, by  taking each sum  we have 
\begin{eqnarray}\nonumber
\sum_{x\in \Sing(\fol)} Log(\fol,\widehat{D}_k, x) &=& \sum_{x\in \Sing(\fol)} Log(\fol,D, x) + \\\nonumber\\\nonumber &+& \left[\sum_{\substack{x\in \Sing(\fol),\\x\in \left(D - \widehat{D}_k \right) \cap D_k }}
 \mu_x(\fol|_{D_k}) +\sum_{\substack{x\in \Sing(\fol),\\x\in \widehat{D}_k \cap D_k}}
 Log(\fol|_{D_k},\widehat{D}_k \cap D_k, x)\right] = \\\nonumber\\\nonumber &=&
\sum_{x\in \Sing(\fol)} Log(\fol,D, x) + \sum_{x\in \Sing(\fol|_{D_k})}  Log(\fol|_{D_k},\widehat{D}_k \cap {D}_k, x),
\end{eqnarray}
\noindent where in the last step we  use  the relation (\ref{eq0_03}) and the following equality of sets $\left(D - \widehat{D}_k \right) \cap D_k = D_k - (\widehat{D}_k \cap D_k)$.
$\square$

\section{Proof Corollary  \ref{cor-Pn}}
By Theorem \ref{teo1_1} we have that 
 $$
  \displaystyle\int_{\mathbb{P}^n}c_n(T_{\PP^n}(-\log\, D)- \mathcal{O}_{\mathbb{P}^n}(1-d))=   \sum_{x\in  \Sing(\fol)\cap (X\setminus  D)} \mu_x(\fol) + \sum_{x\in \Sing(\fol)  \cap  D }    Log(\fol,D,x),
 $$
since $T_{\fol}=\mathcal{O}_{\mathbb{P}^n}(1-d)$, where $d$ is the degree of $\fol$.  We have to prove the positivity of
 $$
  \displaystyle\int_{\mathbb{P}^n}c_n(T_{\PP^n}(-\log\, D)-\mathcal{O}_{\mathbb{P}^n}(1-d)).
 $$

Firstly, we will prove the following formula  
  $$
  \displaystyle\int_{\mathbb{P}^n}c_n(T_{\PP^n}(-\log\, D)-\mathcal{O}_{\mathbb{P}^n}(1-d))=\displaystyle\sum^{n}_{i=0} {n+1\choose i}  \sigma_{n-i}(d_1,\dots,d_k, d-1).
 $$
It follows from  \cite[Theorem 4.3]{EA} that 
\begin{eqnarray} \nonumber
0 \longrightarrow T_{\PP^n}(-\log\, D) \longrightarrow \OO_{\PP^n}(1)^{n+1}\oplus  \OO_{\PP^n}^{k-1} \longrightarrow  \bigoplus_{i=1}^{k} \OO_{\PP^n}(d_i)\longrightarrow 0.
\end{eqnarray}
From this exact sequence we get   
 $$
 c(T_{\PP^n}(-\log\, D))= \frac{(1+h)^{n+1}}{  \prod_{i = 1}^k  (1+d_ih)^i},
 $$
 where $h=c_1(\OO_{\PP^n}(1))$. 
Since $c(T_\fol)=1+(1-d)h$ we have 
 \begin{eqnarray}\nonumber
 c(T_{\PP^n}(-\log\, D)-T_{\fol})&=& \frac{c(T_{\PP^n}(-\log\, D))}{c(T_\fol)} 
 \\\nonumber
 &=& \frac{(1+h)^{n+1}}{ (1+(1-d)h) \prod_{i = 1}^k  (1+d_ih)^i} 
\\\nonumber \\\nonumber &=&  \frac{(1+h)^{n+1}}{  \prod_{i = 0}^{k}  (1+m_ih)^i}, 
\end{eqnarray}
 where $m_i:=d_{i+1}$, for all $i=0,\dots,k-1$ and $m_{k}=1-d$. 
 Therefore
 $$
  c_n(T_{\PP^n}(-\log\, D)-T_{\fol})=\left[\frac{(1+h)^{n+1}}{  \prod_{i = 0}^{k}  (1+m_ih)^i}\right]_n= \displaystyle\sum^{n}_{i=0}(-1)^{n-i} {n+1\choose i}  \sigma_{n-i}(m_0,\dots,m_{k})h^n. 
 $$
Thus, we have that 
  $$
\displaystyle\sum^{n}_{i=0}(-1)^{n-i} {n+1\choose i}  \sigma_{n-i}(m_0,\dots,1-d) =  \displaystyle\sum^{n}_{i=0} {n+1\choose i}  \sigma_{n-i}(m_0,\dots,d-1).   
 $$
Now, suppose that $d=2$  and  $d_i=1$, for all $k$, then 
 $$
c(T_{\PP^n}(-\log\, D)-T_{\fol})=  \frac{(1+h)^{n}}{  (1+h)^{k}}.
 $$ 
Therefore,  the conclusion is the same  as  \cite[Example 1.1]{CSV}. If $d>1$, then it is clear that 
 $$0<\displaystyle\sum^{n}_{i=0} {n+1\choose i}  \sigma_{n-i}(d_1,\dots,d_k, d-1) = \int_{\PP^n} c_n(T_{\PP^n}(-\log\, D)-T_{\fol}).$$  
In order to conclude the item (ii) we recall that $Log(\fol,D,x)=0$ for all non-degenerate singularity $x\in \Sing(\fol) \cap D$. Then
 $$0<\displaystyle\sum^{n}_{i=0} {n+1\choose i}  \sigma_{n-i}(d_1,\dots,d_k, d-1) =\sum_{x\in \singf \cap  D }  Log(\fol,D,x)$$
says us that there exist at least one  degenerate singular point. 
\section{Proof Corollary \ref{singlogres}}

Consider   a  functorial log  resolution  $\pi: (X,D) \to (Y, \emptyset)$ with excepcional divisor $D$, see  \cite[ Theorems 3.35, 3.34]{kollar07}.     
Since the singular locus of $X$ is invariant with respect to any automorphism it follows from   \cite[ Corollary 4.6]{GKK10}) that the vector field $v\in H^0(Y,TY)$ has a lift  
$\tilde{v}\in H^0(X,T_X(-\log D))$. Denoting by  $\tilde{\fol}$ the foliation associated to the vector field $\tilde{v}$, we 
have that $T_{\tilde{\fol}}=\O_X$,  since $\tilde{v} \in H^0(X,T_X(-\log D))\subset H^0(X,T_X)$ is a global vector field.

Now, since $\Sing(v)\cap \Sing(X)=\emptyset $, we have that $\tilde{v}$ has no zeros along the normal crossings divisor  $D$. By Theorem \ref{teo1_1} we have  
  \begin{eqnarray}\nonumber
\displaystyle\int_{X}c_{n}(T_{X}(-\log\, D)- \O_X) =\displaystyle\int_{X}c_{n}(T_{X}(-\log\, D))  &=&    \sum_{\pi^{-1}(x)\in  \Sing(\tilde{v}) \cap (X\setminus  D) } \mu_{\pi^{-1}(x)}(\tilde{v})  \\ \nonumber
\\ \nonumber
&=& \sum_{x\in  \Sing(v)\cap Y_{reg} } \mu_x(v) 
\\ \nonumber
\\ \nonumber
&=&  \sum_{x\in  \Sing(v) } \mu_x(v)
\end{eqnarray}
  \noindent  since  $\mu_{\pi^{-1}(x)}(\tilde{v})=  \mu_x(v)$, for all $x\in  \Sing(v) $.


\begin{thebibliography}{99}


\bibitem{Ale_p} A. G. Aleksandrov, {\it The index of vector fields and logarithmic differential forms}, Funct. Anal. Appl. 39 (4) (2005) 245-255.

\bibitem{Aluffi} P. Aluffi, {\it Chern classes for singular hypersurfaces}, Trans. Am. Math. Soc. 351 (1999), no. 10, 3989-4026.


\bibitem{EA} E. Angelini, {\it Logarithmic bundles of hypersurface arrangements in $\PP^n$}, Collectanea
Mathematica, Volume 65, Issue 3, (2014), 285-302





\bibitem{BB1} P. Baum  and  R. Bott, {\it Singularities of Holomorphic Foliations}, J. Differential Geom, 7 (1972), 279-342.

\bibitem{a03} J.-P. Brasselet, J. Seade and T. Suwa, {\it An explicit cycle representing the Fulton-Johnson class}, Singularit\'es Franco-Japonaises, S\'emin. Congr., 10, Soc. Math. France, Paris, p. 21-38, 2005.



\bibitem{BarSaeSuw} J.-P. Brasselet, J. Seade and   T. Suwa, {\it Vector Fields on Singular Varieties}, Lecture Notes in Mathematics, Spring, 2009.

  \bibitem{Correa-toricas}
M. Corr\^ea Jr,  {\it Darboux--Jouanolou--Ghys integrability for one-dimensional foliations on toric varieties}, Bull. Sci. Math., 134(7):693-704, 2010.



\bibitem{CFLM} M. Corr\^ea,   Antonio M. Ferreira,  F. Louren\c co and   D. Machado, {\it On Gauss-Bonnet and Poincar\'e-Hopf type  theorems for  complex  $\partial$-manifolds}. To appear in Moscow Mathematical Journal, 2020.
(arXiv:1808.05178v1).


\bibitem{C-L}
M. Corr\^ea and  F. Louren\c co, \emph{Determination of Baum-Bott residues for  higher dimensional  foliations}. 
Asian Journal of Mathematics,
vol. 23, n.3, (2019), 527-538.




\bibitem{CoMa} 
M. Corr\^ea and  D. Machado, {\it  GSV-index for holomorphic Pfaff systems}.  Documenta Mathematica. \textbf{25}, 1011-1027, 2020. 


\bibitem{MD} M. Corr\^ea and D. Machado,  {\it Residue formulas for logarithmic foliations and applications}, Trans. Am. Math. Soc. 371, 6403-6420, 2019.

 

\bibitem{CMS1}
M. Corr\^ea Jr., L. G. Maza, and M. G. Soares. {\it Algebraic integrability of polynomial differentialr-forms}.J.Pure Appl. Algebra, 215(9):2290-2294, 2011.

\bibitem{CMS2}
M. Corr\^ea Jr., L. G. Maza, and M. G. Soares. {\it Hypersurfaces invariant by Pfaff systems}.Commun. Contemp.Math., 17(6):1450051, 18, 2015.

\bibitem{CPS}
M. Corr\^ea Jr., M. Rodriguez Pe\~na and M. G. Soares .{\it  A Bott-type residue formula on complex orbifolds}. Int. Math. Res. Not., (10): 2889-2911, 2016.





\bibitem{CSV} F. Cukierman, M.  Soares and   I. Vainsencher, 
{\it  Singularities of logarithmic foliations}. Compositio Math. \textbf{142} 131--142 (2006).

\bibitem{Deligne} P. Deligne, {\it Equations differentielles \`a points singulier r\'eguliers}, Lecture Notes in  Mathematics, 163, Springer-Verlag, 1970.


\bibitem{a08} X. G\'omez-Mont, {\it An algebraic formula for the index of a vector field on a hypersurface
with an isolated singularity}, J. Algebraic Geom. 7 (1998), 731-752.



\bibitem{GKK10}
Daniel Greb, Stefan Kebekus, and S{\'a}ndor~J. Kov{\'a}cs, \emph{Extension
  theorems for differential forms and {B}ogomolov-{S}ommese vanishing on log
  canonical varieties}, Compos. Math. \textbf{146} (2010), no.~1, 193--219.
  \MR{2581247}


\bibitem{GSV}  X. G\'omez-Mont, J. Seade  and  A. Verjovsky, {\it The index of a holomorphic flow with an isolated singularity}, Math. Ann. 291 (1991), 737-751.



 
\bibitem{Lita} S. Iitaka,  {\it Logarithmic forms of algebraic varieties}, J. Fac. Sci. Univ. Tokyo Sect. IA 23 (1976), 525-544.

\bibitem{Katz} N. M. Katz, {\it The regularity theorem in algebraic geometry}, Actes Congres Intern.  Math., 1970, t.1, 437-443.



 

 
 
 \bibitem{kollar07} J{\'a}nos Koll{\'a}r, 
 \emph{Lectures on resolution of singularities}, Annals of Mathematics
  Studies, vol. 166, Princeton University Press, Princeton, NJ, 2007.
  \MR{2289519 (2008f:14026)}
 

 
\bibitem{a07} D. Lehmann, M. Soares  and T. Suwa, {\it On the index of a holomorphic vector field tangent to a singular variety}, Bol. Soc. Bras. Mat. 26 (1995), pp. 183-199.
 
\bibitem{YN} Y. Norimatsu, {\it Kodaira Vanishing Theorem and Chern Classes for $\partial$-Manifolds}, Proc. Japan Acad., 54, Ser. A. (1978), 107-108.

\bibitem{Sai} K. Saito,  {\it Theory of logarithmic differential forms and logarithmic vector fields}, J. Fac. Sci. Univ. Tokyo, 27(2), p. 265-291, 1980.
 
\bibitem{SeaSuw1} J. Seade and T. Suwa, {\it A residue formula for the index of a holomorphic flow}, Math. Ann. 304 (1996), 621-634.
 
\bibitem{RS} R. Silvotti, {\it On a conjecture of Varchenko}, Invent. Math. 126 (1996), no. 2, 235-248.
 

 
\bibitem{Suw2} T. Suwa, {\it Indices of vector fields and residues of singular holomorphic foliations}, Actualit\'es Math\'ematiques, Hermann \'Editeurs des Sciences et des Arts, 1998.

\end{thebibliography}
\end{document}